\documentclass[12pt,twoside,a4paper]{article}
\usepackage[cp1251]{inputenc}
\usepackage{srcltx}
\usepackage[russian]{babel}
\usepackage{amsmath,amsfonts,amscd,amssymb,latexsym}

\begin{document}



\bigskip

\begin{center}

 {\bf Maltsev bases for partially commutative nilpotent groups

\bigskip

 E.I.~Timoshenko}

\bigskip

Department of Algebra and Logic

Novosibirsk State Technical University

20 K. Marx ave, Novosibirsk, Russia, 630073,\\
eitim45@gmail.com

  \bigskip

 Abstract

\end{center}

\begin{quote}

 We construct an ordered set of commutators in a partially
commutative nilpotent group $F(X; \Gamma; \mathfrak N_m)$. This
 set allows us to define a canonical form for each
 element of $F(X; \Gamma; \mathfrak N_m)$. Namely, we
 construct a Maltsev basis for the group $F(X; \Gamma; \mathfrak
 N_m).$

\end{quote}

\bigskip


\begin{center}

{\bf  1. Introduction}

\end{center}

Many classes of algebraic structures are defined through the
category of simple graphs. One of these is the class of so-called
partially commutative algebraic structures.  These structures have
well known applications both in mathematics and in computer sciences
as well as in robotics.

In the paper,  let $\Gamma = \langle X; E\rangle$ be a simple
graph with the set of vertices $X$ and the set of edges $ E.$

First partially commutative structures being studied were monoids.
A free partially commutative monoid  on $X$ associated with
 $\Gamma$ is the monoid denoted by $M(X;
\Gamma)$ which is defined by the monoid presentation
$$M(X; \Gamma)= \langle X; xy =yx , \{x,y\} \in E \rangle.$$
The notion of partially commutative monoid was introduced by P.
Cartier and D. Foata in 1969 (see [1]) to study combinatorial
problems in connection with word rearrangements.

A free partially commutative group $F(X; \Gamma)$ is closely
related to $M(X; \Gamma).$ It is defined by the group presentation
$$F(X; \Gamma) = \langle X; xy=yx, \{x,y\} \in E \rangle.$$

The groups  $F(X; \Gamma)$  were first introduced in the 1970’s by
A. Baudisch (see [2]) as “semifree groups”  and then were studied
in the 1980’s by C. Droms (see [3, 4, 5]) calling these groups by
“graph groups”.

The class of free partially commutative groups contains free and
free abelian groups. Free partially commutative groups possess a
number of remarkable properties. For example, a group
$F(X;\Gamma)$ is a residually torsion-free nilpotent group (see [6]).
Therefore free partially commutative groups are torsion-free.
These groups are
linear (see [7]). Fundamental groups of almost all surfaces are
subgroups of free partially commutative groups (see [8]). In [9],
it is observed that two free partially commutative groups $F(X;
\Gamma)$ and $F(Y; \Delta)$ are isomorphic iff their defining
graphs $\Gamma$ and $\Delta$ are isomorphic.

Free partially commutative  groups have provided several crucial
examples having shaped the theory of finitely presented groups;
notably  M. Bestvina and N. Brady's example of a homologically
finite (of type FP) but not geometrically finite (in fact not of
type F2) group; and Mikhailova's example of a group with
unsolvable subgroup membership problem.  Recently it was shown by
work of M. Sageev, F. Haglund, D. Wise, I. Agol and others that
many well-known families of groups virtually embed into free
partially commutative groups: among these are Coxeter groups,
limit groups, and
fundamental groups of closed 3-manifold groups (see for example
[10]).

Consider a variety $\mathfrak M$ of groups. A partially
commutative group in $\mathfrak M$ with a defining graph $\Gamma$
is a group $F(X;\Gamma; \mathfrak M)$ defined as $$ F(X; \Gamma;
\mathfrak M) = \langle X;  xy=yx, \{x,y\} \in E\rangle$$ in
$\mathfrak M.$

Consequently a free partially commutative group $F(X; \Gamma)$ is
a partially commutative group in the variety $\mathfrak G$ of all
groups.

Recall that the commutator $(g_1,g_2)$ of two elements
$g_1, g_2$ of a group $G$ is defined by $(g_1, g_2) =
g_1^{-1}g_2^{-1}g_1g_2.$  By $\mathfrak A^2$ denote a variety of
all metabelian groups, i.e., all groups satisfying an identity
$((x, y), (z, v)) = 1.$

Among all partially commutative groups
$F(X; \Gamma; \mathfrak M)$, $\mathfrak M \neq \mathfrak G,$ the most studied case is the case of partially commutative metabelian groups, i.e. the groups
$F(X; \Gamma; \mathfrak{A}^2)$.

 Let $v(x_{i_1},\ldots, x_{i_m})$ be a representation of an element $v \in F(X; \Gamma; \mathfrak M)$
 as a product of generators in $X,$ where the vertices $x_{i_1},\ldots, x_{i_m}$ occur in this representation.
Then set $\sigma(v) = \{x_{i_1},\ldots, x_{i_m}\}$.  Denote by
$\Gamma_v$ the subgraph of $\Gamma$  generated by the set
$\sigma(v)$ and by $\Gamma_{v,x}$  the connected component of the
graph $\Gamma_v$ such that this component contains a vertex $x \in
\sigma(v)$.
 Let us order the set $X$ as $x_1 < x_2 < \ldots <
x_n.$ By $\max(\Gamma_{v,x})$ denote the greatest vertex in the
connected component $\Gamma_{v,x}$.

\bigskip

The following theorem describes a basis of the commutant $G'$ of a
partially commutative metabelian group $G = F(X; \Gamma; \mathfrak
A^2)$.

\bigskip

{\bf  Theorem 1.} [11] {\it Let the set $X = \{x_1, \ldots, x_n\}$
of vertices of a graph $\Gamma$ be ordered as follows $x_1 < x_2 <
\ldots < x_n$  and let $G= F(X; \Gamma; \mathfrak A^2).$ Then a
basis of the commutant $G'$ is the set  consisting of all elements
$v$ of the form $$v = u^{-1}(x_i, x_j)u,\,\, \text{where}\,\,\,  u
= x_{j_1}^{t_1}\ldots x_{j_m}^{t_m}, \,\,\,\{t_1,\ldots, t_m\}
\subset \mathbb Z \backslash \{0\},$$ such that the following
conditions hold:

 (a) $j \leq j_1 < j_2 \ldots < j_m \leq n, 1 \leq j < i \leq n;$

 (b) the vertices $x_i$ and $x_j$ are in
different connected components
 of the graph  $\Gamma_v;$

(c) $x_i = \max(\Gamma_{v,x_i}).$}

\bigskip

 There are results obtained for centralizers and annihilators  of groups $F(X; \Gamma; \mathfrak A^2)$ ([12]), embeddings these groups into matrix groups (see [13]), and their groups of automorphisms ([14]) and values of centralizer dimensions ([15, 16]). The universal and
elementary theories of these groups are investigated in [12, 17,
18].

The lower central series of a group $G$ is the sequence of
subgroups $G_{(n)},\ n \geq 1,$ defined inductively as follows
$$G_{(1)} =
G, \,\,\, G_{(i+1)} = (G_{(i)}, G),$$ where $(G_{(i)}, G)$ denotes
the subgroup of $G$ generated by the commutators $(x, y) $ with $x
\in G_{(i)}, y \in G.$

A variety $\mathfrak N_c$ consists of all groups $G$ such that
$G_{(c+1)} = 1.$

The properties of partially commutative nilpotent groups $F(X;
\Gamma; \mathfrak N_c)$  are much less studied.  Even a canonical
form for elements of groups $F(X; \Gamma; \mathfrak N_c)$ for $c
\geq 4$ is not known yet (the cases $c = 2, 3$ are considered in
[19]).

In this paper, we study partially commutative groups in $\mathfrak N_c.$ All groups considered below are finitely generated. So, the set $X=\{x_1,\ldots,x_n\}$ is finite.

For a subset $H$ of a group $G$ denote by $gp\langle H\rangle$ the subgroup generated by $H$.

 If $G$ is a torsion-free finitely generated
nilpotent group then $G$ has a central series  $$G = G_1 > G_2 >
\ldots > G_{s+1} = 1 \eqno (1)$$ with infinite cyclic factors.
Take elements $a_1,\ldots,a_s$ such that $G_i = gp\langle
a_i,G_{i+1}\rangle.$

\bigskip

{\bf Definition 1.} (see [20]) {\it An ordered system
$\{a_1,\ldots, a_s\}$ of elements is called a Maltsev basis for
$G$ obtained by  the  central series (1).}

\bigskip

The construction of a Maltsev basis of a group makes it possible
to indicate a canonical form of its elements. Every element $g \in
G$ can be  uniquely represented in the form $$g = a_1^{t_1} \ldots
a_s^{t_s},\,\, t_i \in \mathbb Z.$$

A Maltsev basis for a group $F(X; \Gamma; \mathfrak A^2 \wedge
\mathfrak N_c)$ was found in [19]. Let us recall its description.

Define a commutator $c_m =
(y_1, y_2,\ldots,y_m)$, where
$y_i \in X,$ by induction: $c_2 = (y_1,
y_2), c_m = (c_{m-1},y_m).$

Let $B$ be the set of commutators of the form $$v = (x_{j_1},
x_{j_2}, \ldots, x_{j_m}),\,\, 2 \leq m \leq c,$$ in a
group $F(X; \Gamma; \mathfrak A^2 \wedge \mathfrak N_c)$ such that the following conditions hold:

(a) $1 \leq j_2 \leq  j_3 \leq  j_m \leq  n,\,\, j_2 < j_1\leq n;$

(b) the vertices $x_{j_1}$  and $x_{j_2}$  are in different
connected components of the graph $\Gamma_v;$

 (c) $x_{j_1}  = \max\,(\Gamma_{v,x_{j_1}} ).$

 \bigskip

  {\bf Theorem 2.} [19] {\it The
set of elements $ X \sqcup B$ is a Maltsev basis of  $F(X; \Gamma,
\mathfrak A^2 \wedge \mathfrak N_c)$  obtained  by refining the
lower central series of this group. }

\bigskip

The group $F(X; \Gamma; \mathfrak N_c) \cong
F(X;\Gamma)/F_{(c+1)}(X; \Gamma)$ is torsion-free (see [6],
Theorem 2.1). This means that there exists a Maltsev basis for
$F(X; \Gamma; \mathfrak N_c).$

\bigskip

The aim of this paper  is to find a Maltsev basis for the group
$F(X;  \Gamma; \mathfrak N_c).$

The study of the free partially commutative Lie algebra was
started by G. Duchamp in 1987 (see [21]). Let $R$ be a  domain. A
free partially commutative Lie $R$-algebra $\mathcal
L_R(X;\Gamma)$ is the Lie algebra defined by the Lie algebra
presentation  $$\mathcal
L_R(X;\Gamma) = \langle X; [x_i, x_j]= 0, \{x_i, x_j\} \in
E\rangle.$$

Put $\mathcal L(X;\Gamma) = \mathcal L_{\mathbb Z}(X;\Gamma).$ In [6], the relation between
the graded Lie
$\mathbb Z$-algebra associated with the quotients of the lower
central series of $F(X; \Gamma)$ and the Lie algebra
$\mathcal L(X; \Gamma)$ was established. We are going to use this relation.

 The concept of basic commutators was
introduced by Ph. Hall in [22]. Hall's commutators are usually
used in group theory.

For convenience, we will use so called standard commutators (see
[23]) for the description of a Maltsev basis.

  Denote by $X^*$  the set of all   words in
$X = \{x_1,\ldots, x_n\}$ including the empty word denoted by 1.
We also denote by $|u|$ the length of any  $u \in X^*.$ Let us
extend an arbitrary linear order  on $X$ to a lexicographic order
''<'' on $X^*$  as follows.  Put $u < 1$  for each $ 1 \neq u \in
X^*$ and by induction put $x_iu' < x_jv'$ if $x_i < x_j$ or $x_i =
x_j, u' < v'$.

\bigskip

{\bf Definition 2.} {\it Let $$ALS(X) = \{ u \in
X^*\,\,|\,\,\forall u_1, u_2 \in X^*(u = u_1u_2 \longrightarrow
u_2u_1 <u_1u_2)\}.$$ A word $u \in ALS(X)$ is called an associative
Lyndon---Shirshov word.}

\bigskip

Let us define a set $ \mathcal G(X)$ and a bar map $ \mathcal
G(X) \longrightarrow X^*$ as follows.

\bigskip

{\bf Definition 3.} {\it (a) $x_i \in \mathcal G(X)$ for all $x_i
\in X,\,\,\,\overline{x_i} = x_i.$

(b) If $u, v \in \mathcal G(X),$ then $(u, v) \in \mathcal G(X)$
and $\overline{(u,v)} = \overline{u}\,\,\overline{v}.$}

\bigskip

The bar map erases all parentheses and commas.

We put $$\mathcal G_m(X)  = \{u\,\,|\,\,u \in \mathcal
G(X),\,\,|\overline{u}| = m\}, \text{ then }
\mathcal G(X) = \bigcup_{m\geq 1}\mathcal G_m(X).$$

Now we give a definition of the set $(X^*)$ of standard
commutators.

\bigskip

{\bf Definition 4.} {\it (a) $x_i \in (X^*)$ for $i=1,\ldots,n.$

(b) Let $w = (u, v).$ Then $w \in (X^*)$    if and only if the
following  conditions are true:

(b1) $\overline{w} \in ALS(X);$

(b2) $u, v \in (X^*),\, \overline{u} > \overline{v} ;$

(b2) if $u= (u_1, u_2)$ then $\overline{v} \geq \overline{u_2}.$}

\bigskip

Let $$(X^*)_m = \{u \,|\,u \in (X^*), |\overline{u}|= m\}.$$

If $F$ is the free  group with the basis $X=\{x_1,\ldots,x_n\},$
and $(x, y) = x^{-1}y^{-1}xy$ for $x, y \in F,$ then the set of
commutators $(X^*)_m$ forms a basis of the free abelian group
$F_{(m)}/F_{(m+1)}$ for $m=1,2 \ldots$ (see [23], Theorem 3.5).

\bigskip

{\bf Definition 5.} {\it Let $u \in X^*$. By $\delta_i(u)$ denote the number of occurrences of $x_i$ in $u.$ For $u \in X^*,$
put
$$supp(u) = \{x_i \,|\,\delta_i(u) \neq 0\}.$$}

Finally, let us define by induction a subset $\mathcal C(X;
\Gamma)$ of $(X^*).$

\bigskip

{\bf Definition 6.} {\it (a) All elements of $X$  belong to
$\mathcal C(X; \Gamma).$

(b) An element $u \in (X^*)_m, m \geq 2,$ belongs to $\mathcal
C(X, \Gamma)$ if $u = (v,w),$ where $v$ and $w$ are elements of
$\mathcal C(X; \Gamma)$ and there is an element in $supp(v)$ such
that this element is not connected in $\Gamma$ with the first
letter of $w.$

(c) There are no other elements in $\mathcal C(X; \Gamma).$}

\bigskip

Let $$\mathcal C_i(X; \Gamma) = \{u \in \mathcal C(X; \Gamma)\,|\,
|\overline{u}| = i,\,\,\, i = 0,1,\ldots\}.$$


Let ''$\prec$'' on $\mathcal C(X; \Gamma)$ such that $ u \prec v$
if $u \in \mathcal C_p(X; \Gamma), \,\, v \in \mathcal C_q(X;
\Gamma),\,\,1 \leq p <q.$

Let $$\mathcal C^{(m)}(X; \Gamma) = \bigcup_{1\leq i \leq m}
\mathcal C_i(X; \Gamma).$$

\bigskip

{\bf Theorem 3.} {\it The set  $\mathcal C^{(m)}(X; \Gamma)$ with
respect to the order ''$\prec$'' is a Maltsev basis for  the group
$F(X; \Gamma; \mathfrak N_m)$ obtained  by refining the lower
central series.}

\bigskip

\begin{center}

 {\bf 2. Bases for partially commutative nilpotent Lie algebras}

 \end{center}

An explicit construction for bases of free partially commutative
Lie algebras was obtained in [24]. To give this description let us
first recall a definition of Lyndon---Shirshov words.

The lexicographic order ''<''  has been defined above (in the last paragraph before
 Definition 2as well as the set $ALS(X)$ of associative Lyndon---Shirshov words on $X$ (see Definition 2).

Let $\mathcal L(X)$ be the free Lie algebra on the set  $X=
\{x_1,\ldots,x_n\}.$

Let us give a definition of a set $[X^*],\,\, [X^*] \subseteq
\mathcal L(X),$ of non-associative Lyndon---Shirshov words.


\bigskip

{\bf Definition 7.} {\it
 (a) $x_i \in [X^*]$ for $i=1,\ldots, n;$

 (b)  Let $[w] = [[u], [v]].$  Then $[w] \in [X^*]$ if and only if the following
 conditions hold:

 (b1) $w \in ALS(X)$;

 (b2) $[u], [v] \in [X^*],\, u> v,$ where $u, v$ denote the words
 in $X^*$ obtained from $[u], [v]$ by omitting the Lie brackets $[\,,\,];$

 (b3) if $[u] = [[u_1], [u_2]]$ then $v \geq u_2.$}

\bigskip

 It was shown in [25]   that the set $[X^*]$ of all  non-associative  Lyndon---Shirshov words  is a linear  basis of the
free Lie $R$-algebra $\mathcal L_R(X)$ over a domain
$R$.

For a free partially commutative Lie algebra $\mathcal L_R(X;
\Gamma)$ over a domain $R$ define inductively the set of partially
commutative Lyndon---Shirshow words (PCLS-words for short)  by induction.

\bigskip

{\bf Definition 8.} {\it (a) All elements of $X$ are PCLS-words.

(b) A Lyndon-Shirshov word $[u]$ such that $|u| >1$ is a PCLS-word
if $[u] = [[v], [w]],$ where $[v]$ and $[w]$ are PCLS-words and
there is an element in $supp(v)$ such that it is not connected in
$\Gamma$ with the first letter of $w.$

(c) There are no other PCLS-words.}

 \bigskip

 Denote the set of all PCLS-words of a free partially
commutative Lie $R$-algebra $\mathcal L_R(X; \Gamma)$ by $PCLS(X;
\Gamma).$

The first  result on bases of free partially commutative Lie
algebras $\mathcal L_R(X; \Gamma)$ was obtained by D. Duchamp and
D. Krob in [26], but they did not give an explicit description of
a basis.

Using the method of Gr\"{o}bner---Shirshov bases E. Poroshenko in
[24] obtained an explicit description of bases for free partially
commutative Lie algebras.

\bigskip

 {\bf Theorem 4.} [24] {\it  Let $R$ be a unital commutative ring and
 $\Gamma$ be a graph. Then the set $PCLS(X;\Gamma)$ is a linear basis of the free partially
commutative Lie $R$-algebra $\mathcal L_R(X;\Gamma).$}

\bigskip

Let $\mathfrak L_m$ be a variety of  all nilpotent Lie algebras of
 class at most $m.$ Denote by $\mathcal L_R(X; \Gamma; \mathfrak
L_m)$ the partially commutative $m-$nilpotent $R$-algebra Lie.

A linear basis for a partially commutative nilpotent Lie algebra
can be easily obtained from a linear basis for the corresponding
free partially commutative Lie algebra.

\bigskip

{\bf Theorem 5.} [24] {\it Let $R$ be a unital commutative ring and
  $\Gamma$ be a graph.
 Then a linear basis of the partially commutative nilpotent
$R$-algebra $\mathcal L_R(X; \Gamma;\mathfrak L_m)$  consists of
all elements of $PCLS(X;\Gamma)$ whose lengthes are not greater
than $m.$}

\bigskip


\begin{center}

 {\bf 3. Proof of Theorem 3}

 \end{center}

Let $G$ be a group. Define the associated graded abelian group
 $gr(G)$ as follows
 $$gr(G) =
\bigoplus_{m \geq 1}gr_m(G),$$ where $gr_m(G) = G_{(m)}/G_{(m
+1)},$. The group $gr(G)$ has a structure of a graded Lie algebra
over the ring $\mathbb Z$ of integers with the bracket operation
in $gr(G)$ induced by the commutator operation in $G$ .







By $\mathcal F$ denote a graded $\mathbb Z$-module
$$\mathcal F = \bigoplus_{m \geq 1} F_{(m)}(X;\Gamma)/F_{(m+1)}(X;\Gamma).$$


The element $g^* \in \mathcal{F}$ is called a homogeneous element of degree $m$ if
this element is in $F_{(m)}(X;\Gamma)/F_{(m+1)}(X;\Gamma)$.

$\mathcal F$ can be equipped with a Lie $\mathbb Z$-algebra structure as follows.

Let $g^*$ and $h^*$ be homogeneous elements of  degrees $m$ and
$n$ respectively. Denote by $g$ a preimage of $g^*$ in
$F_{(m)}(X;\Gamma)$ and by $h$ a preimage of $h^*$ in
$F_{(n)}(X;\Gamma).$ Then  $(g, h) \in F_{(m+n)}(X;\Gamma)$
according to a property of the lower central series. Thus we can
equip $\mathcal F$ with the Lie bracket defined by the relation
$$[g^*, h^*] = (g, h)F_{(m+n+1)}(X;\Gamma).$$ This Lie bracket is well-defined.
 It does not depend on the choice of preimages $g$ and $h$ for elements $g^*$ and
 $h^*.$ We can  extend the bracket operation to
 $\mathcal F$ by distributivity.

Let vertices $x_i$ and $x_j$ be adjacent in $\Gamma.$ Then
$$[x_i F_{(2)}(X;\Gamma), x_j F_{(2)}(X;\Gamma)] = (x_i, x_j)F_{(3)}(X;\Gamma)
= F_{(3)}(X;\Gamma) =0$$ in $\mathcal F.$ Therefore,  we can
extend mapping
$$\alpha(x) = x F_{(2)}(X;\Gamma),\,\,x \in X,$$ to a homomorphism
of the Lie $\mathbb Z$-algebras $\mathcal L(X; \Gamma)$ and
$\mathcal F$:
 $$\alpha : \mathcal L(X; \Gamma)
\longrightarrow \mathcal F.$$

Let us now define a family $\mathcal A(X)$ of
$\mathcal L(X, \Gamma)$ by induction. We set  $\mathcal A_1(X) = X.$ For  $m
\geq 2,$ put
$$\mathcal A_m(X) = \{[u,v]\,|\,u \in  \mathcal A_p(X),\,v \in \mathcal
A_q(X),\,\,p+q =m\},$$ $$\mathcal A(X) = \bigcup_{m \geq 1}
\mathcal A_m(X). $$

Let $\mathcal L_m(X; \Gamma)$ be a submodule of $\mathcal L(X;
\Gamma)$ generated by $\mathcal A_m(X).$

In  [6], Theorem 2.1, it was proved that $\alpha$ is an isomorphism of
graded Lie algebras from $\mathcal L(X; \Gamma)$ graded by
$(\mathcal L_m(X; \Gamma)_{m \geq 1})$ into $\mathcal F.$

Consequently
 $$F_{(m)}(X; \Gamma)/F_{(m+1)}(X; \Gamma) \stackrel
{\alpha}{\simeq} \mathcal L_m(X; \Gamma),$$ for $m \geq 1.$

By $PCLS_m(X; \Gamma)$ denote the set of all  $PCLS(X; \Gamma)$
words of  length $m.$ As follows from Theorem 5, the set
$PCLS_m(X; \Gamma)$ forms a basis for the additive abelian group
$\mathcal L_m(X; \Gamma).$

Comparing Definitions 4 and 7, and then Definition 6 and 8, we see
that the isomorphism $\alpha$  maps the set $PCLS_m(X; \Gamma)$
onto the set $\mathcal C_m(X; \Gamma).$ Therefore, $\mathcal
C_m(X; \Gamma)$ forms a basis for the abelian group $F_{(m)}(X;
\Gamma)/F_{(m+1)}(X; \Gamma).$
This completes the proof.

\bigskip

{\bf Example.} Let $\Gamma = \langle x_1, x_2, x_3; \,\,\{x_1,
x_2\}\rangle ,$  $x_1 > x_2 > x_3.$

By construction,
 $$\mathcal C^{(3)}(X; \Gamma) = \{x_1, x_2, x_3;  (x_1, x_3), (x_2, x_3);   (x_1,(x_1, x_3)),$$$$  (x_2,(x_2,x_3)), ((x_1,x_3),x_2),
((x_1,x_3),x_3), ((x_2,x_3),x_3)\}$$ is a Maltsev basis of the
group $F(X;\Gamma; \mathfrak N_3).$

\bigskip

{\bf References}

\bigskip

[1]\, P. Cartier, D. Foata, Problemes combinatoires de commutation
et de rearrangements,  Lecture Notes in Mathematics {\bf
85}(1969), Springer-Verlag, Berlin, New York.

[2]\, A. Baudisch, Subgroups of semifree groups, Acta Math. Acad.
Sci. Hungar {\bf 38}(1-4)(1981) 19–28.

[3]\, C. Droms, Graph groups, coherence, and three-manifolds, J.
Algebra, {\bf 106}(2) (1987) 484–489.

[4]\, C. Droms, Isomorphisms of graph groups, Proc. Amer. Math.
Soc. {\bf 100}(3) (1987) 407–408.

[5]\, C. Droms, Subgroups of graph groups, J. Algebra {\bf 110}(2)
(1987) 519–522.

[6]\, G. Duchamp, D. Krob, The lower central series of the free
partially commutative group, Semigroup Forum {\bf 45}(1992)
385-494.

[7]\, T. Hsu and D. Wise, On linear and residual properties of
graph products, Mich. Math. J. {\bf 46}(2)(1999) 251-259.

 [8]\, J. Crisp and B. Wiest, Embeddings of graph braid groups and surface groups
in right-angled Artin groups and braid groups, Alg., Geom., Topol.
{\bf 4}(2004) 439-472.

[9]\, C. Droms, Isomorphisms of graph groups, Proc. Am. Math. Soc.
{\bf 100}(1987) 407-408.

[10]\,  D.T. Wise, From riches to raags: 3-manifolds, right-angled
Artin groups, and cubical geometry, CBMS Regional Conference
Series in Mathematics, vol. 117, Published for the Conference
Board of the Mathematical Sciences, Washington, DC; by the
American Mathematical Society, Providence, RI, 2012. MR 2986461.

[11]\, E. I. Timoshenko, A basis of partially commutative
metabelian groups, Izv. Math. {\bf 85} (2021) (to appear).

[12]\, Ch.K. Gupta, E.I.Timoshenko, Partially commutative
metabelian groups: centralizers and elementary equivalence,
Algebra and Logic{\bf 48}(3) (2009), 173-192.

[13]\, E. I. Timoshenko, On embedding of partially commutative
metabelian groups to matrix groups, Inter. J. of Group Theory {\bf
7}(4) (2018), 17-26.

[14]\, E. I. Timoshenko, Automorphisms of partially commutative
metabelian groups, Algebra and Logic {\bf 59}(2) (2020)165-179.

[15]\, E. I. Timoshenko, Centralizers dimensions and universal
theories for partially commutative metabelian groups, Algebra and
Logic {\bf 56}(2) (2017), 149-170.

[16]\, E. I. Timoshenko, Centralizers dimensions of partially
commutative metabelian groups, Algebra and Logic {\bf 57}(1)
(2018) 69-80.

[17]\, Ch. K. Gupta, E.I. Timoshenko, Universal theories for
partially commutative metabelian groups, Algebra and Logic {\bf
50}(1) (2011) 1-16.

[18]\, E. I. Timoshenko, Universal equivalence of partially
commutative metabelian groups, Algebra and Logic {\bf 49}(2)
(2010) 177-196.

[19]\, E. I. Timoshenko, A Maltsev basis for a partially
commutative nilpotent metabelian group, Algebra and Logic {\bf
50}(5) (2011) 439-446.

[20]\, V. I. Kargapolov, Yu. I. Merzlyakov, Foundations of the
group theory,  Moscow (Russian).

[21]\, G. Duchamp, Algorithms sur les polin\^{o}mes en variables
non commutatives, Th\`{e}se d'Universit\`{e}, Universit\`{e} Paris
7, LITP Report No.87-58  (1987) 87-58.

[22]\, Ph. Hall, Some word problems, J. London Math. Soc., {\bf
33} (1958), 482-496.

[23]\, K.T. Chen, R. H. Fox, R. C. Lyndon,  Free Differential
Calculus, IV. The Quotient Groups of the Lower Central Series.
 The Annals of Mathematics, 2nd Ser. {\bf 68}(1) (1958) 81-95.

 [24]\, E. N. Poroshenko, Bases for partially commutative Lie
algebras, Algebra and Logic {\bf 50}(5) (2011) 405-417.

[25]\,  A. I. Shirshov, On Free Rings, Math. Sb. {\bf 45(87)}(2)
(1958) 113-122 (Russian).

 [26]\, G. Duchamp, D. Krob, The Free Partially Commutative Lie
Algebra: Bases and Ranks, Advances in Mathematics {\bf 92} (1992)
95-126.

\end{document}